\documentclass[10pt]{article}
\usepackage{graphicx}
\usepackage{amsmath,amssymb,amsthm,amsfonts}
\usepackage{amssymb}
\usepackage{mathtools}
\newtheorem{thm}{Theorem}[section]

\newtheorem{lemma}[thm]{Lemma}

\theoremstyle{definition}

\theoremstyle{remark}
\newtheorem{rem}{Remark}[section]

\numberwithin{equation}{section}

\DeclareMathSymbol{\C}{\mathalpha}{AMSb}{"43}

\newcommand{\bsub}{\begin{subequations}}
\newcommand{\esub}{\end{subequations}$\!$}

\begin{document}

\title{A variational problem arising in registration of diffusion tensor images \thanks{Email: hanhuan11@mails.ucas.ac.cn (H. Han); hszhou@wipm.ac.cn (H. Zhou)
}\thanks{This work was supported by NSFC under grant No.11471331 and partially supported by National Center for Mathematics and Interdisciplinary Sciences.}}

\author{Huan Han$^a$
 \ and\  Huan-Song Zhou$^b$\\
\small $^a$Wuhan Insitute of Physics and Mathematics, Chinese Academy of Sciences,\\ \small  P.O. Box 71010, Wuhan, 430071, China
\\
\small$^b$Department of Mathematics, Wuhan University of Technology, Wuhan 430070, China}

\date{}

\smallbreak \maketitle

\begin{abstract} The existence of a global minimizer for a variational problem arising in registration of diffusion tensor images is proved, which ensures that there is a regular spatial transformation for the registration of diffusion tensor images.
\end{abstract}

\vskip 0.2truein



\vskip 0.2truein

\section{Introduction}
Let $\Omega\subset \mathbb{R}^3$ be an open bounded domain with Liptschitz boundary $\partial \Omega$, and let $T$ and $D$ be two functions such that
\begin{align}\label{intro1}
T:\Omega\rightarrow \mathbb{R}^d\ \ \ \ \mathrm{and } \ \ \ \ D:\Omega\rightarrow \mathbb{R}^d\ \ \ (d\geq 1).
\end{align}

By talking about the image registration, $T$ and $D$ are viewed as two images in a spatial domain $\Omega$, one is called the floating image(e.g. $T$) and the other is the target image (e.g. $D$). The registration of two images is to find a smooth and locally non-degenerate spatial transformation $h:\Omega\rightarrow \Omega$ such that the composition of $T$ and $h$, that is $T\circ h(x)=T(h(x))$, is close to $D$ in some sense, for example, find a $h$ such that the sum of squared errors in a suitable space reaches to its minimum. Note that the dimension $d$ in (\ref{intro1}) of the space for the ranges of $T$ and $D$ varies with imaging modalities, for examples, images acquired by Magnetic Resonance Imaging(MRI) or Computed Tomography(CT) are scalar valued with $d=1$, RGB images are vector valued with $d=3$, and the Diffusion Tensor Imaging(DTI) is matrix valued with $d=9$ (or $d=6$ if the matrix at each voxel of DTI is  $3\times 3$ symmetric positive).

In order to get a spatial transformation $h:\Omega\rightarrow \Omega$ with higher regularity, Dupuis, Grenander and Miller\cite{ref10} improved the variational model proposed by Amit\cite{ref4} for scalar image registration ($d=1$), and they considered the following variational problem in a suitable space for $v$:
\begin{align}\label{intro2}
&\hat{v}=\arg \mathop {\min }\limits_{v} \int_0^\tau\|Lv(\cdot,t)\|^2_{L^2(\Omega)}dt+\|T\circ h(\cdot)-D(\cdot)\|^2_{L^2(\Omega)},&
\end{align}
where $\tau>0$ is the time duration, ``$\circ$" is the composition of two functions, and $L:[H_0^{3}(\Omega)]^3\rightarrow [L^2(\Omega)]^3$ is a differential operator satisfying
\begin{align}\label{intro3}
&\|Lv(\cdot,t)\|^2_{L^2(\Omega)}=\sum_{i=1}^3\int_{\Omega}|(Lv)_i(x,t)|^2dx\geq c \|v(\cdot,t)\|^2_{[H^3_0(\Omega)]^3}, &
\end{align}
for some $c>0$ and for all $ v\in [H^3_0(\Omega)]^3 \triangleq H^3_0(\Omega) \times H^3_0(\Omega) \times H^3_0(\Omega)$ with $H^3_0(\Omega)=W^{3,2}_0(\Omega)$. Clearly, $L=\sum\limits_{i+j+k=3}\frac{\partial^3}{\partial x_1^i\partial x_2^j\partial x_3^k}$ satisfies (\ref{intro3}).  Furthermore, $v$ and $h$ are constrained by the following equation:
\begin{align} \label{intro4}
\frac{d{\eta}(s;t,x)}{ds}\triangleq\dot{\eta}(s;t,x)=v(\eta(s;t,x),s),\ \eta(t;t,x)=x, \   h(x)=\eta(0;\tau,x), \ 0\leq s,t\leq \tau,&
\end{align}
here $\eta(s;t,x)$ means that a particle placed at $y=\eta(s;t,x)$ at time $s$ is transformed to a point $x$ at time $t$ under the forcing term $v(x,t)$. In \cite{ref10}, the authors gave a rigid mathematical proof on the existence of global minimizer, which validates the applications of model (\ref{intro2}) in numerical simulations of scalar image registration. Based on \cite{ref10}, the large derivation principle(LDP) of the constraint equation (\ref{intro4}) are concerned in \cite{new}. However, the model (\ref{intro2}) does not work for registration of DTI images because each voxel of DTI image contains a $3\times 3$ symmetric positive definite real matrix (i.e., diffusion tensor) and the orientation of diffusion tensors must be considered in making spatial transformation, which is much more complicated than that of scalar images. For this reason, there are two tensor reorientation strategies, the finite strain(FS) and the preservation principle direction(PPD), have been proposed by Alexander\cite{ref3}, which are widely used for analyzing DTI data. Note that, for DTI images, both $T$ and $D$ are maps from $\Omega$ to the set of $3\times 3$ Symmetric Positive Definite real matrixes ($SPD(3)$ in short), that is,
\begin{align}\label{intro5}
T,D:\Omega\rightarrow SPD(3)\subset \mathbb{R}^6.
\end{align}

In order to extend the variational model (\ref{intro2}) to the case of registration of DTI images, based on FS reorientation strategy, Li et al. \cite{ref151} introduced a new transformation operation ``$\diamond$" which is used to replace the usual composition operation ``$\circ$" in (\ref{intro2}). The operation ``$\diamond$"  is given by
\begin{align}\label{intro6}
T\diamond h(x)=R \left[T\circ h(x)\right] R^T \ \ \mathrm{with} \ \ R=J_x^T\left(J_xJ_x^T\right)^{-\frac{1}{2}}\mathrm{and} \ \ J_x=\nabla _x h^{-1}(x),
\end{align}
and the computation of $\left(J_xJ_x^T\right)^{-\frac{1}{2}}$ is given in Appendix.

With the above definition and notations, the variational model proposed in \cite{ref151} for the registration of DTI images can be formulated by
\begin{align}\label{intro7}
&\hat{v}=\arg \mathop {\min }\limits_{v \in \mathcal{F}} H(v),&
\end{align}
where $H(v)=\int_0^\tau\|Lv(\cdot,t)\|^2_{L^2(\Omega)}dt+\|T\diamond h(\cdot)-D(\cdot)\|^2_{L^2(\Omega)}$, $\tau$ and $L$ are the same as in (\ref{intro2}), $h(x)$ is given by (\ref{intro4}), $\mathcal{F}$ is defined by (\ref{intro11}).

\begin{rem} The constrain variational problem(\ref{intro7}) and (\ref{intro4}) is well defined, see Lemma \ref{th-5}. Moreover,
$J_x$ defined in (\ref{intro6}) satisfies
\begin{align}\label{intro8}
&J_x=\nabla _x h^{-1}(x)=\nabla _x \eta(\tau;0,x).&
\end{align}
Note that here function $h^{-1}$ is well defined, see Remark \ref{newremark}.

In fact, let $y=h(x)=\eta(0;\tau,x)$, then by definition of $\eta(s;t,x)$, we know $x=h^{-1}(y)=\eta(\tau;0,y)$. Hence, for any particle at position $y$ of the moving image at time $0$, the position at time $\tau$ of this particle is $x$. By [1, Section III], we know $J_y=\nabla_y x$, that is,
\begin{align}\label{eq-p}
&J_y=\nabla_y \eta(\tau;0,y).   &
\end{align}\qed
\end{rem}

The authors of \cite{ref151} designed an efficient numerically algorithm for registration of DTI images based on (\ref{intro7}) with $L\equiv 0$. If $L\neq 0$, the situation becomes much more complicated. To ensure the model (\ref{intro7}) can be well used in developing an optimal algorithm of DTI registration, it is important to know  mathematically the solvability of (\ref{intro7}). To the authors' knowledge, there seems  no any results on the existence of global minimizer for variational problem (\ref{intro7}). Motivated by \cite{ref10}, the aim of this paper is to give a rigid mathematical proof on the existence.

Before giving our main result, we introduce some notations.

For a given $\tau>0$ and $u(x,t):\Omega\times[0,\tau]\rightarrow \mathbb{R}^3$, let $L$ be a linear differential operator given by (\ref{intro3}), we define
\begin{align}\label{intro11}
\mathcal{F}\triangleq \{u(x,t):u(x,t)\in [H^3_0(\Omega)]^3\  \mathrm{for\ any} \ t\in[0,\tau] \ \mathrm{and} \int_0^\tau\|Lu(\cdot,t)\|^2_{L^2(\Omega)}dt<+\infty\},
\end{align}
endowing with the following inner product and norm
\begin{align}\label{intro12}
(u,v)_{{\mathcal{F}}}=\int_0^\tau \int_\Omega Lu(x,t)\cdot Lv(x,t)dxdt \ \mathrm{and}\ \|u\|_{\mathcal{F}}^2=\int_0^\tau\|Lu(\cdot,t)\|^2_{L^2(\Omega)}dt,
\end{align}
and $\mathcal{F}$ is a separable Hilbert space associated with the above norm.

Throughout the paper, the norm of a matrix $A(x)=\left(
                                                   \begin{array}{c}
                                                     a_{ij}(x) \\
                                                   \end{array}
                                                 \right)_{n\times m}
$ is defined by
\begin{align}\label{intro13}
&\|A(x)\|=\sum_{i=1}^{n}\sum_{j=1}^{m} |a_{ij}(x)|. &
\end{align}

We say a matrix $A(x)$ is continuous at $x\in \Omega$ if each element $a_{ij}(x)$ of $A(x)$ is continuous at $x$.
Particularly, if $m=1$, the matrix becomes an $n\times 1$ vector which is usually denoted $v(x)=(v_1(x),v_2(x),\cdots,v_n(x))^T$ with the norm $\|v(x)\|=\sum\limits_{i=1}^n |v_i(x)|$.

Based on the above definitions and notations, our main result can be stated as follows:

\begin{thm}\label{th-A10}
For $\Omega$ being given in (\ref{intro1}), let $T$ and $D$ be the two maps defined by (\ref{intro5}), and let the set $\triangle_T \triangleq \{x\in \Omega: T(\cdot) \ \ \mathrm{ is \ discontinuous\ at} \ x\}$ be a set of measure zero. If $L:[H_0^{3}(\Omega)]^3\rightarrow [L^2(\Omega)]^3$ is a linear differential operator and satisfies (\ref{intro3}) and
$$\mathop {\max }\limits_{x\in \Omega} \|T(x)\|<+\infty,\ J=\mathop {\max }\limits_{x\in \Omega} \|T(x)-D(x)\|^2<+\infty,$$
 then the variational problem (\ref{intro7}) has a global minimizer $\hat{v}\in {\mathcal{F}}$, which induces a deformation $\hat{h}(x)\in [C^{1,\frac{1}{2}}(\Omega)]^3$ from $\Omega$ to $\Omega$. Moreover, the derivative of $\hat{h}(x)$ satisfies (\ref{eq-5}).
\end{thm}

\section{Preliminary results}
In this section, we show some lemmas which are required in proving Theorem \ref{th-A10}.
\begin{lemma}\label{lem-2}
For any $f=(f_1, f_2, f_3)^T\in [H^3_0(\Omega)]^3$, there exists $K>0$ such that
\vskip1mm {\bf (i)}  $\|f(x)-f(y)\|\leq K\|Lf\|_{L^2(\Omega)}\|x-y\|$, $\forall x,y\in \Omega$.
\vskip1mm {\bf(ii)}  $\|\nabla f(x)-\nabla f(y)\|\leq K\|Lf\|_{L^2(\Omega)}\|x-y\|^{\frac{1}{2}}$,  $\forall x,y\in \Omega$, where $\nabla f=\left(
                  \begin{array}{c}
                    f_{ij} \\
                  \end{array}
                \right)_{3\times 3}
$ with $f_{ij}=\frac{\partial f_i}{\partial x_j}$ $(i,j=1,2,3)$.
\end{lemma}
\noindent{\bf Proof:}\ {\bf (i)} By Taylor's formula, there holds
\begin{align}
&\|f(x)-f(y)\|=C\|\nabla f(\zeta)\cdot(x-y)^T\|,\nonumber&
\end{align}

 Note that $H^3(\Omega)\hookrightarrow C^{1,\frac{1}{2}}(\Omega)$, then
 \begin{align}
&\|f(x)-f(y)\|\leq C\|\nabla f(\zeta)\|\|x-y\|\leq \tilde{C}\|f\|_{C^{1,\frac{1}{2}}(\Omega)}\|x-y\|\leq K\|Lf\|_{L^2(\Omega)}\|x-y\|.\nonumber&
\end{align}

{\bf (ii)}  Note that $H^{3}(\Omega)\hookrightarrow C^{1,\frac{1}{2}}(\Omega)$. Then,
\begin{align}
&\frac{\|\nabla f_i(x)-\nabla f_i(y)\|}{\|x-y\|^{\frac{1}{2}}}\leq \|f_i\|_{C^{1,\frac{1}{2}}(\Omega)}\leq C\|f_i\|_{H^{3}_0(\Omega)}\leq \tilde{K}\|Lf\|_{L^2(\Omega)}.\nonumber&
\end{align}
\vskip1mm That implies,
\begin{align}
&\|\nabla f_i(x)-\nabla f_i(y)\|\leq  \tilde{K}\|Lf\|_{L^2(\Omega)}\|x-y\|^{\frac{1}{2}}. \nonumber&
\end{align}
\vskip1mm Hence, we obtain that
\begin{align}
&\|\nabla f(x)-\nabla f(y)\|=\sum_{i=1}^3\|\nabla f_i(x)-\nabla f_i(y)\|\leq  K\|Lf\|_{L^2(\Omega)}\|x-y\|^{\frac{1}{2}}. \ \ \ \ \ \ \ \ \ \ \ \ \ \ \ \ \ \ \ \ \ \ \ \ \ \ \ \ \ \ \qed\nonumber&
\end{align}

\vskip1mm By (\ref{intro7}), we know that $H(v)$ is essentially a functional of $v$ and $\eta$, where $v$ and $\eta$ satisfy (\ref{intro4}). If we write $H(v)$  as a functional about $v$, then the existence and uniqueness for the solutions of (\ref{intro4}) should be known. Otherwise, the definition of (\ref{intro7}) will be ambiguous if one $v$ produces two or more $\eta$ by (\ref{intro4}).
So, the following results on the ODE problem (\ref{intro4}) imply that the variational problem (\ref{intro7}) and (\ref{intro4}) is well defined.
\begin{lemma}\label{th-5}
Let $v\in {\mathcal{F}}$ and $\|v\|^2_{\mathcal{F}}<+\infty$. If $v(\cdot,t)|_{\mathbb{R}^3\setminus{\Omega}}=0$ for each $t\in[0,\tau]$, then there exists  a unique solution $\eta(s;t,x)\in C([0,\tau],\bar{\Omega})$ of (\ref{intro4}) for $s\in [0,\tau]$ and $x\in \bar{\Omega}$.
\end{lemma}
\noindent{\bf Proof:}
 First, we choose a constant $M$ large enough such that $\|v\|^2_{\mathcal{F}}\leq M$. For $x\in \overline{\Omega}$ and $s\in[0,\tau]$, we define
\begin{align}\label{bgh}
& \Gamma_s(\phi)=x+\int_t^s v(\phi(r),r)dr \ \ \mathrm{for}\  \phi\in C([0,\tau],\mathbb{R}^{3})&
\end{align}

Then, we have

 \begin{align}\label{nere}
  \|\Gamma_s(\phi)\|\leq& \|x\|+\left\|\int_t^s v(\phi(r),r)dr\right\|\leq \|x\|+\left|\int_0^\tau \|v(\phi(r),r)\|dr\right|\nonumber&\\
  &\leq \|x\|+\left|\int_0^\tau \|v(\cdot,r)\|_{C(\Omega)}dr\right|\leq \|x\|+c_1\left|\int_0^\tau \|v(\cdot,r)\|_{H^3(\Omega)}dr\right|\nonumber&\\
  &\leq \|x\|+c\left|\int_0^\tau \|Lv(\cdot,r)\|_{L^2(\Omega)}dr\right|\leq \|x\|+c\tau^{\frac{1}{2}}\left(\int_0^\tau \|Lv(\cdot,r)\|^2_{L^2(\Omega)}dr\right)^{\frac{1}{2}}\nonumber&\\
  &\leq E+c(M\tau)^{\frac{1}{2}}\triangleq \bar{M},&
\end{align}
where $E=\sup\{\|x\|:x\in\Omega\}$.

Now we claim that the mapping $\Gamma_s$ from $C([0,\tau],\mathbb{R}^{3})$ to $C([0,\tau],\mathbb{R}^{3})$ is continuous. In fact, by the definition of  $\Gamma_s$ and Lemma \ref{lem-2}, we obtain that
\begin{align}
&\|\Gamma_s(\phi_1)-\Gamma_s(\phi_2)\|
\leq \left|\int_t^s \|v(\phi_1(r),r)-v(\phi_2(r),r)\|dr\right|  \nonumber&\\
\leq& K\left|\int_t^s \|Lv(\cdot,r)\|_{L^2(\Omega)}\|\phi_1(r)-\phi_2(r)\| dr\right|
\leq K\left|\int_t^s \|Lv(\cdot,r)\|_{L^2(\Omega)} dr\right|\|\phi_1-\phi_2\|_{C([0,\tau],\overline{\Omega})}  \nonumber&\\
\leq &K\tau^{\frac{1}{2}}\left(\int_0^\tau \|Lv(\cdot,r)\|^2_{L^2(\Omega)} dr\right)^{\frac{1}{2}}||\phi_1-\phi_2||_{C([0,\tau],\overline{\Omega})}
\leq K(M\tau)^{\frac{1}{2}}\|\phi_1-\phi_2\|_{C([0,\tau],\overline{\Omega})}.  \nonumber&
\end{align}

 Hence,
\begin{align}
& \|\Gamma_s(\phi_1)-\Gamma_s(\phi_2)\|_{C([0,\tau],\overline{\Omega})}\leq K(M\tau)^{\frac{1}{2}}||\phi_1-\phi_2||_{C([0,\tau],\overline{\Omega})},  \nonumber&
\end{align}
and  the claim is proved.

Next, we show that $\Gamma_s$ from $C([0,\tau],\mathbb{R}^{3})$ to $C([0,\tau],\mathbb{R}^{3})$ is compact.
In fact, for $0\leq t_1<t_2\leq \tau$, we see that
\begin{align}
\|\Gamma_{t_1}(\phi)-\Gamma_{t_2}(\phi)\|&=\left\|\int_{t_1}^{t_2}v(\phi(r),r) dr\right\|
\leq |t_2-t_1|^{\frac{1}{2}}\left(\int_{t_1}^{t_2}\|v(\phi(r),r)\|^2 dr\right)^{\frac{1}{2}} \nonumber&\\
&\leq |t_2-t_1|^{\frac{1}{2}}\left(\int_{t_1}^{t_2}\|v(\cdot,r)\|_{C^{1,\frac{1}{2}}(\Omega)}^2 dr\right)^{\frac{1}{2}}
\leq C|t_2-t_1|^{\frac{1}{2}}\left(\int_{t_1}^{t_2}\|v(\cdot,r)\|_{H_0^{3}(\Omega)}^2 dr\right)^{\frac{1}{2}} \nonumber&\\
&\leq K|t_2-t_1|^{\frac{1}{2}}\left(\int_{t_1}^{t_2}\|Lv(\cdot,r)\|_{L^2(\Omega)}^2 dr\right)^{\frac{1}{2}}
\leq K|t_2-t_1|^{\frac{1}{2}}\left(\int_{0}^{\tau}\|Lv(\cdot,r)\|_{L^2(\Omega)}^2 dr\right)^{\frac{1}{2}} \nonumber&\\
&\leq KM^{\frac{1}{2}}|t_2-t_1|^{\frac{1}{2}}, \nonumber&
\end{align}
this means that $\{\Gamma_s(\phi),\phi\in C([0,\tau],\mathbb{R}^{3})\}$ is equicontinuous for $s\in [0,\tau]$.
Then the Arzela-Ascoli Theorem\cite{ref161} shows that $\{\Gamma(\phi):\phi\in C([0,\tau],\mathbb{R}^{3})\}$ is relative compact, which means that any bounded sequence has as convergent subsequence. Hence, the mapping $\Gamma_s:C([0,\tau],\mathbb{R}^{3})\rightarrow C([0,\tau],\mathbb{R}^{3})$ is compact.
On the other hand, the set
\begin{align}\label{newread}
\Lambda \triangleq \{\phi\in C([0,\tau],\mathbb{R}^{3}): \phi=\lambda\Gamma_s({\phi})  \ \ \mathrm{for\ \ some}\ 0\leq \lambda\leq 1\}
\end{align}
is not empty and bounded. In fact, it follows from (\ref{bgh}) that
\begin{align}\label{ngkl}
\|\lambda\Gamma_s(\phi_1)-\lambda\Gamma_s(\phi_2)\|=&\lambda\left\|\int_t^s[v(\phi_1(r),r)-v(\phi_1(r),r)]dr\right\|\leq\lambda\left|\int_0^\tau K\|Lv(\cdot,r)\|_{L^2(\Omega)}\|\phi_1(r)-\phi_2(r)\|dr\right|\nonumber&\\
\leq&\lambda K\int_0^\tau \|Lv(\cdot,r)\|_{L^2(\Omega)}dr\|\phi_1(\cdot)-\phi_2(\cdot)\|_{C([0,\tau],\mathbb{R}^{3})}\nonumber&\\
\leq&\lambda K\tau^{\frac{1}{2}}\left(\int_0^\tau \|Lv(\cdot,r)\|^2_{L^2(\Omega)}dr\right)^{\frac{1}{2}}\|\phi_1(\cdot)-\phi_2(\cdot)\|_{C([0,\tau],\mathbb{R}^{3})}\nonumber&\\
\leq&\lambda K(M\tau)^{\frac{1}{2}}\|\phi_1(\cdot)-\phi_2(\cdot)\|_{C([0,\tau],\mathbb{R}^{3})},&
\end{align}
this implies that $\lambda\Gamma_s$ is a strict contraction
if $\lambda\in(0,\frac{1}{K(M\tau)^{\frac{1}{2}}})$, and the Banach's fixed point theorem shows that $\Lambda\not= \emptyset$.
 Furthermore, for any $\tilde{\phi}\in C([0,\tau],\mathbb{R}^{3})$ satisfies
\begin{align}
\tilde{\phi}=\lambda\Gamma_s(\tilde{\phi}),  0\leq \lambda\leq 1,
\end{align}
it follows from (\ref{nere}) that $\|\tilde{\phi}\|=\|\lambda\Gamma_s(\tilde{\phi})\|\leq \lambda\|\Gamma_s(\tilde{\phi})\| \leq \widetilde{M}$, that is, the set $\Lambda$
is bounded, and the claim is proved.

With the above facts, the Schaefer's Fixed Point Theorem[6, Theorem 4 in Section 9.2] shows  that $\Gamma_s$ has a fixed point $\phi(s)$ such that
\begin{align}
&\phi(s)=\Gamma_s(\phi)=x+\int_t^s v(\phi(r),r)dr\in C([0,\tau],\bar{\Omega}). \nonumber&
\end{align}

 Define $\eta(s;t,x)=\phi(s)$, then we know (\ref{intro4}) has a solution.

 Finally, we prove the uniqueness of the solution of (\ref{intro4}).
\vskip1mm Assume there $\eta_1(s;t,x),\eta_2(s;t,x)$ are two different solutions of (\ref{intro4}), then
\begin{align}
\|\eta_1(s;t,x)-\eta_2(s;t,x)\|&=\left\|\int_t^s v(\eta_1(r;t,x),r)-v(\eta_2(r;t,x),r)dr \right\|\nonumber&\\
&\leq \left|\int_t^s \|v(\eta_1(r;t,x),r)-v(\eta_2(r;t,x),r)\|dr \right|\nonumber&\\
&\leq \left|\int_t^s K\|Lv\|_{L^2(\Omega)}\|\eta_1(r;t,x)-\eta_2(r;t,x)\|dr \right|.\nonumber&
\end{align}
\vskip1mm By Grownwall inequality[10, Lemma 1.1], we know that
\begin{align}
\|\eta_1(s;t,x)-\eta_2(s;t,x)\|&= 0.\nonumber&
\end{align}
\vskip1mm Hence, $\eta_1(s;t,x)\equiv\eta_2(s;t,x)$.\qed

\vskip1mm Our following lemma shows that $J_x$ is not singular.
\begin{lemma}\label{th-jk6}
Let $v\in {\mathcal{F}}$ and the condition (\ref{intro3}) be satisfied. Then
\vskip1mm {\bf (i)} For any $t\in [0,\tau],s\in [0,\tau]$ and $x\in \bar{\Omega}$, $\eta(s;t,x)$ is differentiable with respect to (w.r.t. in short) $x$ and the derivative $\nabla_x\eta(s;t,x)\triangleq \Theta(s;t,x)$ satisfies
 \begin{equation}\label{eq-5}\left\{\begin{aligned}
 &\dot{\Theta}(s;t,x)=\nabla_\eta v(\eta(s;t,x),s)\Theta (s;t,x), &\\
 &\Theta (t;t,x)=I=(\delta_{ij})_{3\times 3}, &\\
 \end{aligned}\right.\end{equation}
 where $\delta_{ij}=0$ if $i\neq j$, $\delta_{ij}=1$ if $i= j$,  and\  $\nabla_\eta v(\eta(s;t,x),s)=\left(
                                                                                  \begin{array}{c}
                                                                                    \frac{\partial v_i}{\partial \eta_j} (\eta(s;t,x),s)\\
                                                                                  \end{array}
                                                                                \right)_{3\times 3}
 $.
\vskip1mm  {\bf(ii)} The determinant of $\Theta(s;t,x)$ is equal to
\begin{align}
&{\rm det}(\Theta(s;t,x))=e^{\int_t^s\sum\limits_{i=1}^3 v_{i,\eta_i}(\eta(s;t,x),s)ds}. \nonumber&
\end{align}
\end{lemma}
\noindent{\bf Proof:} {\bf (i)} By (\ref{intro3}) and the fact that $H^3_0(\Omega)\hookrightarrow C^{1,\frac{1}{2}}(\Omega)$, we know
\begin{align} \label{ode}
&\dot{\eta}(s;t,x)=v(\eta(s;t,x),s),\ \ \ \ \eta(t;t,x)=x,  &
\end{align}
and $v(\eta(s;t,x),s)$ is differentiable w.r.t $\eta(s;t,x)$ and $\eta(s;t,x)$ is differentiable w.r.t. $x$ (cf. \cite{ref16}).
\vskip1mm Differentiating  the first equation of (\ref{ode}) w.r.t. $x$, we have
\begin{align}
\frac{d}{ds}{\nabla_x\eta}(s;t,x)&=\nabla_x(v(\eta(s;t,x),s))
 =\nabla_\eta v(\eta(s;t,x),s)\cdot\nabla_x\eta(s;t,x),\nonumber
\end{align}
that is,
\begin{align}\label{eq-6}
&\dot{\Theta}(s;t,x)=\nabla_\eta v(\eta(s;t,x),s)\cdot\Theta(s;t,x).&
\end{align}
If we differentiate  the second equation of (\ref{ode}) w.r.t. $x$, which gives
\begin{align}\label{eq-7}
&\Theta(t;t,x)=I,&
\end{align}
that is, the second equation of (\ref{eq-5}). So, part {\bf (i)} is proved.

\vskip1mm {\bf (ii)} Let $\theta_{ij},v_{ij}$ be the element at the $i th$ row and the $j th$ column of $\Theta(s;t,x)$, $\nabla_\eta v(\eta(s;t,x),s)$ respectively, then we obtain that
\begin{align}
\frac{d \det(\Theta(s;t,x))}{ds}&=(v_{11}+v_{22}+v_{33})\det(\Theta(s;t,x)).\nonumber&
\end{align}
\vskip1mm Hence,
\begin{align}
&{\rm det}(\Theta(s;t,x))={\rm det}(I)e^{\int_t^sv_{11}+v_{22}+v_{33}ds}=e^{\int_t^s\sum\limits_{i=1}^3 v_{i,\eta_i}(\eta(s;t,x),s)ds}.\ \ \ \ \ \ \ \ \ \ \ \ \ \ \ \ \ \ \ \ \ \ \ \ \ \ \ \ \ \ \ \ \ \ \qed\nonumber&
\end{align}
\begin{rem}\label{newremark}
Let $t=\tau$ and $s=0$, by Lemma \ref{th-jk6}{\bf (ii)}, there holds
\begin{align}\label{ensuer}
\det(\nabla h(x))=\det(\Theta(0;\tau,x))=e^{-\int_0^\tau\sum\limits_{i=1}^3 v_{i,\eta_i}(\eta(s;t,x),s)ds}\neq 0,
\end{align}
for all $x\in \Omega$. By Inverse Function Theorem[5, Theorem 7 in Appendix], (\ref{ensuer}) ensures the existence of $h^{-1}$.
\end{rem}

\begin{lemma}\label{new}
For $\Theta_{n}(n=1,2,3,\cdots),\Theta\in SPD(3)$, let $\lambda_{n}^{(1)}\geq \lambda_{n}^{(2)}\geq \lambda_{n}^{(3)}>0$ and $\lambda^{(1)}\geq \lambda^{(2)}\geq \lambda^{(3)}>0$ be eigenvalues of $\Theta_{n}$ and $\Theta$, respectively. If $\Theta_{n}\xrightarrow{n} \Theta$ (i.e. each element of $\Theta_{n}$ converges to that of $\Theta$), then $\lambda_n^{(i)}\xrightarrow{n} \lambda^{(i)}(i=1,2,3)$.
\end{lemma}
\noindent{\bf Proof:} Define $\Theta_n=\left(
                                     \begin{array}{c}
                                       \theta^n_{ij} \\
                                     \end{array}
                                   \right)_{3\times3}
$,
$\Theta=\left(
                                     \begin{array}{c}
                                       \theta_{ij} \\
                                     \end{array}
                                   \right)_{3\times3}
$.
Let
\begin{align} \label{eq-Finnew}
& {\rm det}(\tilde{\lambda} I -\Theta_{n})=0\ \ \mathrm{and}\ \ {\rm det}(\lambda I -\Theta)=0.&
\end{align}
\vskip1mm Since $\Theta_{n},\Theta\in SPD(3)$, we know that the equations of (\ref{eq-Finnew}) have the following form
\begin{align}\label{eq-Finnew1}
& \tilde{\lambda}^3+p_{n}\tilde{\lambda}^2+q_{n}\tilde{\lambda}+r_{n}=0,{\lambda}^3+{p}{\lambda}^2+{q}{\lambda}+{r}=0, &
\end{align}
where $p_{n}=-(\theta^{n}_{11}+\theta^{n}_{22}+\theta^{n}_{33}),q_{n}=\theta^{n}_{11}\theta^{n}_{22}-\theta^{n}_{12}\theta^{n}_{21}+\theta^{n}_{11}\theta^{n}_{33}-\theta^{n}_{13}\theta^{n}_{31}+\theta^{n}_{22}\theta^{n}_{33}-\theta^{n}_{23}\theta^{n}_{32}$,
$r_{n}=-\theta^{n}_{11}\theta^{n}_{22}\theta^{n}_{33}+\theta^{n}_{11}\theta^{n}_{23}\theta^{n}_{32}+\theta^{n}_{12}\theta^{n}_{21}\theta^{n}_{33}-\theta^{n}_{12}\theta^{n}_{23}\theta^{n}_{31}-\theta^{n}_{13}\theta^{n}_{21}\theta^{n}_{32}+\theta^{n}_{13}\theta^{n}_{22}\theta^{n}_{31}$.
$p,q,r$ are of the same form as $p_{n},q_{n},r_{n}$, but $\theta^{n}_{ij}$ is replaced by $\theta_{ij}$. Therefore,
\begin{align}
p_n\xrightarrow{n} p, q_n\xrightarrow{n} q, r_n\xrightarrow{n} r.
\end{align}

\vskip1mm Since $\lambda_{n}^{(1)}\geq \lambda_{n}^{(2)}\geq \lambda_{n}^{(3)}>0$ and $\lambda^{(1)}\geq \lambda^{(2)}\geq \lambda^{(3)}>0$ are eigenvalues of $\Theta_{n}$ and $\Theta$ respectively, it follows from (\ref{eq-Finnew1}) that
\begin{align} \label{eq-Finnew3}
& (\tilde{\lambda}-\lambda_n^{(1)})(\tilde{\lambda}-\lambda_n^{(2)})(\tilde{\lambda}-\lambda_n^{(3)})=0,\ \ (\lambda -\lambda^{(1)})(\lambda -\lambda^{(2)})(\lambda -\lambda^{(3)})=0,&
\end{align}
and $p_{n}=-(\lambda_n^{(1)}+\lambda_n^{(2)}+\lambda_n^{(3)})$, $q_{n}=\lambda_n^{(1)}\lambda_n^{(2)}+\lambda_n^{(1)}\lambda_n^{(3)}+\lambda_n^{(2)}\lambda_n^{(3)}$, $r_{n}=-\lambda_n^{(1)}\lambda_n^{(2)}\lambda_n^{(3)}$, $p=-(\lambda^{(1)}+\lambda^{(2)}+\lambda^{(3)})$, $q=\lambda^{(1)}\lambda^{(2)}+\lambda^{(1)}\lambda^{(3)}+\lambda^{(2)}\lambda^{(3)}$, $r=-\lambda^{(1)}\lambda^{(2)}\lambda^{(3)}$.
\vskip1mm Furthermore, for $i=1,2,3$, using (\ref{eq-Finnew1}) we see that
\begin{align} \label{eq-Finnew4}
&(\lambda^{(i)}-\lambda_n^{(1)})(\lambda^{(i)}-\lambda_n^{(2)})(\lambda^{(i)}-\lambda_n^{(3)})=(\lambda^{(i)})^3+p_n(\lambda^{(i)})^2+q_n\lambda^{(i)}+r_n&\nonumber\\
=&[(\lambda^{(i)})^3+p(\lambda^{(i)})^2+q\lambda^{(i)}+r]+(p_n-p)(\lambda^{(i)})^2+(q_n-q)\lambda^{(i)}+r_n-r&\nonumber\\
=&(p_n-p)(\lambda^{(i)})^2+(q_n-q)\lambda^{(i)}+r_n-r\xrightarrow{n} 0.&
\end{align}
\vskip1mm Then, either $\lambda^{(i)}-\lambda_n^{(1)}\xrightarrow{n} 0$, or $\lambda^{(i)}-\lambda_n^{(2)}\xrightarrow{n} 0$, or $\lambda^{(i)}-\lambda_n^{(3)}\xrightarrow{n} 0$. Note that $\lambda_n^{(1)}+\lambda_n^{(2)}+\lambda_n^{(3)}=p_n\xrightarrow{n} p=\lambda^{(1)}+\lambda^{(2)}+\lambda^{(3)}$, $\lambda_n^{(1)}\lambda_n^{(2)}+\lambda_n^{(1)}\lambda_n^{(3)}+\lambda_n^{(2)}\lambda_n^{(3)}=q_n\xrightarrow{n} q=\lambda^{(1)}\lambda^{(2)}+\lambda^{(1)}\lambda^{(3)}+\lambda^{(2)}\lambda^{(3)}$ and $\lambda_n^{(1)}\lambda_n^{(2)}\lambda_n^{(3)}=-r_n\xrightarrow{n} -r=\lambda^{(1)}\lambda^{(2)}\lambda^{(3)}$.
Therefore, $\lambda_n^{(i)}\xrightarrow{n} \lambda^{(i)}(i=1,2,3)$.\qed
\begin{lemma}\label{newtg}
Let $C_m$, $B_m$ be $3\times 3$ matrixes with $\left|\det(B_m)\right|\geq M_0> 0$ and $\|B_m\|\leq M_1<+\infty$. If $C_mB_m\xrightarrow{m} O_{3\times 3}$ , then $C_m\xrightarrow{m} O_{3\times 3}$, where $O_{3\times 3}$ denotes the $3\times 3$ zero matrix (i.e. all elements of the matrix are $0$).
\end{lemma}
\noindent{\bf Proof:} Let $C_m=\left(
                                  \begin{array}{c}
                                    c^m_{ij} \\
                                  \end{array}
                                \right)_{3\times 3}
$ and $B_m=\left(
                                  \begin{array}{c}
                                    b^m_{ij} \\
                                  \end{array}
                                \right)_{3\times 3}$.
Since $C_mB_m\xrightarrow{m} O_{3\times 3}$, $B_m^TC_m^T\xrightarrow{m} O_{3\times 3}$, that is,
\begin{align}
B^{T}_mC_m^i\xrightarrow{m} O_{3\times 1}      \ \ \ \ \ i=1,2,3,
\end{align}
where $C_m^i=(c^m_{i1},c^m_{i2},c^m_{i3})^T$.
\vskip1mm Therefore, there exists $\varepsilon^m_{i1}\xrightarrow{m} 0,\varepsilon^m_{i2}\xrightarrow{m} 0 ,\varepsilon^m_{i3}\xrightarrow{m} 0$  such that
\begin{align}
B^{T}_m C_m^i=(\varepsilon^m_{i1},\varepsilon^m_{i2},\varepsilon^m_{i3})^T      \ \ \ \ \ i=1,2,3.
\end{align}
\vskip1mm  By Cramer's Rule[11, Section 5.3], we obtain that
\begin{align}\label{newh1}
c^m_{i1}=\frac{\varepsilon^m_{i1}B_{m}^{11*}+\varepsilon^m_{i2}B_{m}^{12*}+\varepsilon^m_{i3}B_{m}^{13*}}{\det(B_m)},
\end{align}
\begin{align}\label{newh2}
c^m_{i2}=\frac{\varepsilon^m_{i1}B_{m}^{21*}+\varepsilon^m_{i2}B_{m}^{22*}+\varepsilon^m_{i3}B_{m}^{23*}}{\det(B_m)},
\end{align}
\begin{align}\label{newh3}
c^m_{i3}=\frac{\varepsilon^m_{i1}B_{m}^{31*}+\varepsilon^m_{i2}B_{m}^{32*}+\varepsilon^m_{i3}B_{m}^{33*}}{\det(B_m)},
\end{align}
where $B_{m}^{ij*}$ is the cofactor of $ith$ row and $jth$ column element of $B_m$.
\vskip1mm Since $|B_{m}^{ij*}|\leq 2\|B_m\|^2\leq 2M_1$, it follows from (\ref{newh1}),(\ref{newh2}) and (\ref{newh3}) that
\begin{align}
|c^m_{ij}|\leq&\frac{2\|B_m\|^2(|\varepsilon^m_{i1}|+|\varepsilon^m_{i2}|+|\varepsilon^m_{i3}|)}{|\det(B_m)|}
\leq \frac{2M^2_1(|\varepsilon^m_{i1}|+|\varepsilon^m_{i2}|+|\varepsilon^m_{i3}|)}{M_0}\xrightarrow{m} 0\ \ \ \ i,j=1,2,3.&
\end{align}
\vskip1mm Therefore, $C_m\xrightarrow{m} O_{3\times 3}$.\qed

\begin{lemma}\label{pd7}
Let $\{v_n\in \mathbb{R}^3\}$ be any sequence of ${\mathcal{F}}$ and
\begin{align}\label{eq-d9}
&\|v_n\|^2_{\mathcal{F}}=\int_0^\tau \|Lv_n(\cdot,t)\|^2_{L^2(\Omega)}dt\leq M <+\infty,&
\end{align}
and  $v_n(\cdot,t)|_{\partial \Omega}=0$ for each $t\in[0,\tau]$, then
\vskip1mm {\bf (i)}There exists
a subsequence $\{v_{n_k}\}$ of $\{v_n\}$ and some $v\in {\mathcal{F}}$ such that $v_{n_k}\xrightharpoonup{k} v$ weakly in ${\mathcal{F}}$. Moreover,
\begin{align}\label{eq-d10}
&M\geq\lim_{n_k\rightarrow \infty} \inf \int_0^\tau \|Lv_{n_k}(\cdot,t)\|^2_{L^2(\Omega)}dt\geq \int_0^\tau \|Lv(\cdot,t)\|^2_{L^2(\Omega)}dt, &
\end{align}
\vskip1mm {\bf (ii)}Let $\{v_{n_k}\}$ be the subsequence obtained in {\bf (i)} and let $v$ be its weak limit. For any given $t\in[0,\tau]$, consider the following equations,
 \begin{align}\label{eq-d11}
 \dot{\eta}_{n_k}(s;t,x)=v_{n_k}(\eta_{n_k}(s;t,x),s) \  \ \mathrm{with}\ \ \eta_{n_k}(t;t,x)=x,
 \end{align}
 and
 \begin{align}\label{eq-d12}
 \dot{\eta}(s;t,x)=v(\eta(s;t,x),s) \  \ \mathrm{with}\ \ \eta(t;t,x)=x,
 \end{align}
 then, these two equations have  unique solution in $C([0,\tau]:\bar{\Omega})$ for all $x\in \Omega$. Moreover, for each $s\in [0,\tau]$, $\eta_{n_k}(s;t,x)\in C^{1,\frac{1}{2}}(\Omega)$, $\eta(s;t,x)\in C^{1,\frac{1}{2}}(\Omega)$ with $\eta_{n_k}(s;t,x)\xrightarrow{k} {\eta}(s;t,x)$ uniformly on $[0,\tau]$, and
 \begin{align}\label{eq-d13}
&h_{n_k}(x)=\eta_{n_k}(0;\tau,x)\xrightarrow{k} \eta(0;\tau,x)=h(x) \ \ {\rm for} \ \ \forall x \in \Omega .&
\end{align}
\vskip1mm {\bf (iii)} Let $\Theta_{n_k}(s;0,x)=\nabla_x\eta_{n_k}(s;0,x)$ and $\Theta(s;0,x)=\nabla_x\eta(s;0,x)$, where at $t=0$ $\eta_{n_k}(s;0,x)$ and $\eta(s;0,x)$ are solutions of (\ref{eq-d11}) and (\ref{eq-d12}), respectively. Then $\Theta_{n_k}(s;0,x)\xrightarrow{k} \Theta(s;0,x)$ in the norm (\ref{intro13}) for all $(x,s)\in \Omega\times [0,\tau]$.
\vskip1mm {\bf (iv)} Let $R_{n_k}=\Theta_{n_k}(\Theta_{n_k}\Theta_{n_k}^T)^{-\frac{1}{2}}$ and $R=\Theta(\Theta\Theta^T)^{-\frac{1}{2}}$, here we simply denote $\Theta_{n_k}(s;0,x)$ and $\Theta(s;0,x)$ by $\Theta_{n_k}$ and $\Theta$, respectively. Then $R_{n_k}\xrightarrow{k} R$ in the norm (\ref{intro13}) for all $(x,s)\in \Omega\times [0,\tau]$.
\end{lemma}
\noindent{\bf Proof:} {\bf (i)} and {\bf (ii)} can be proved by the ideas of \cite{ref10} and the Sobolev inequalities[6, Theorem 6 in Section 5.6]. Here, we focus on the proofs of (iii) and (iv).
\vskip1mm {\bf (iii)} By Lemma \ref{th-jk6}, we know that $\Theta_{n_k}(s;0,x)$ and $\Theta(s;0,x)$ are the solution of
\begin{align}\label{eq-newj5k}
\dot{\Theta}_{n_k}(s;0,x)=\nabla_{\eta_{n_k}} v_{n_k}(\eta_{n_k}(s;0,x),s)\Theta_{n_k} (s;0,x),
\ \ \mathrm{and}\ \ \Theta_{n_k} (0;0,x)=I,
\end{align}
 and
 \begin{align}\label{eq-negw5}
 \dot{\Theta}(s;0,x)=\nabla_\eta v(\eta(s;0,x),s)\Theta (s;0,x),\ \ \mathrm{and}\ \ \Theta (0;0,x)=I,
 \end{align}
 respectively.

First, we claim that there exists  $\tilde{M}$ such that
\begin{align}\label{forget}
\|\Theta_{n_k}(s;0,x)\|\leq \tilde{M}<+\infty, \|\Theta(s;0,x)\|\leq \tilde{M}<+\infty.
\end{align}

 In fact, by (\ref{eq-newj5k}),
\begin{align}
\Theta_{n_k}(s;0,x)=I+\int_0^s \nabla_{\eta_{n_k}} v_{n_k}(\eta_{n_k}(s;0,x),s) \Theta_{n_k}(s;0,x)ds.  \nonumber
\end{align}

Hence,
\begin{align}
\|\Theta_{n_k}(s;0,x)\|\leq\|I\|+\int_0^s \|\nabla_{\eta_{n_k}} v_{n_k}(\eta_{n_k}(s;0,x),s)\| \|\Theta_{n_k}(s;0,x)\|ds.  \nonumber
\end{align}
\vskip1mm By the Gronwall inequality[10, Lemma 1.1], we know that
\begin{align}\label{cvb3}
\|\Theta_{n_k}(s;0,x)\|&\leq\|I\|e^{\int_0^s \|\nabla_{\eta_{n_k}} v_{n_k}(\eta_{n_k}(s;0,x),s)\| ds } \leq\|I\|e^{K\int_0^\tau \|Lv_{n_k}(\cdot,s)\|_{L^2(\Omega)} ds } \nonumber& \\
&\leq\|I\|e^{K\tau^{\frac{1}{2}}\left(\int_0^\tau \|Lv_{n_k}(\cdot,s)\|^2_{L^2(\Omega)} ds\right)^{\frac{1}{2}} } \leq\|I\|e^{K\tau^{\frac{1}{2}}M^{\frac{1}{2}} }\doteq \tilde{M}. &
\end{align}

Similarly, we have $\|\Theta(s;0,x)\|\leq \tilde{M}$. This leads to the (\ref{forget}).

Next, we claim that
\begin{align}\label{newadd}
\nabla_x\left[\int_0^t [v_{n_k}(x,r)-v(x,r)]dr\right]\xrightarrow{k} 0\ \ \mathrm{uniformly \ in\ } x\in\Omega\times[0,\tau].
\end{align}

Let $w_{n_k}(x,t)\triangleq v_{n_k}(x,t)-v(x,t)$, $z_{n_k}(x,t)\triangleq\int_0^t w_{n_k}(x,t) dt $. Since $v_{n_k}(x,t)\xrightharpoonup{k} v(x,t)$,  $w_{n_k}(x,t)\xrightharpoonup{k} 0$ and $\|w_{n_k}\|_{\mathcal{F}}\leq 2M$ for all $n_k$.
\vskip1mm If $x,y\in \Omega$ and $s,t\in [0,\tau]$, then
\begin{align}
 &\|z_{n_k}(x,s)-z_{n_k}(y,t)\|\leq\left\|\int_0^s [w_{n_k}(x,r)-w_{n_k}(y,r)]dr\right\|+\left\|\int_s^t w_{n_k}(y,r)dr\right\| \nonumber &\\
 =&\left\|\int_0^s \{[v_{n_k}(x,r)-v_{n_k}(y,r)]-[v(x,r)-v(y,r)]\}dr\right\|+\left\|\int_s^t w_{n_k}(y,r)dr\right\| \nonumber &\\
 \leq &K\left\|\int_0^\tau \|Lv_{n_k}(\cdot,r)\|_{L^2(\Omega)}dr\right\|\|x-y\|
 +K\left\|\int_0^\tau \|Lv(\cdot,r)\|_{L^2(\Omega)}dr\right\|\|x-y\|+\left\|\int_s^t w_{n_k}(y,r)dr\right\| \nonumber &\\
 \leq &K\tau^{\frac{1}{2}}\left(\int_0^\tau \|Lv_{n_k}(\cdot,r)\|_{L^2(\Omega)}^2dr\right)^{\frac{1}{2}}\|x-y\|
 +K\tau^{\frac{1}{2}}\left(\int_0^\tau \|Lv(\cdot,r)\|_{L^2(\Omega)}^2dr\right)^{\frac{1}{2}}\|x-y\| +\left\|\int_s^t w_{n_k}(y,r)dr\right\| \nonumber & \\
 \leq& 2KM^{\frac{1}{2}}\tau^{\frac{1}{2}}\|x-y\|+\left\|\int_s^t w_{n_k}(y,r)dr\right\|
  \leq 2KM^{\frac{1}{2}}\tau^{\frac{1}{2}}\|x-y\|+|t-s|^{\frac{1}{2}}\|w_{n_k}(\cdot,r)\|_{\mathcal{F}} \nonumber & \\
  \leq &2KM^{\frac{1}{2}}\tau^{\frac{1}{2}}\|x-y\|+2M|t-s|^{\frac{1}{2}}. \nonumber &
\end{align}

Hence, $\{z_{n_k}\}\subset C(\Omega\times [0,\tau])^3 $ is equicontinuous, and the Arzela-Ascoli Theorem\cite{ref16} implies that $\{z_{n_k}\}$ is relative compact in $C(\Omega\times [0,\tau])^3$, that is, any bounded sequence has a convergent subsequence in $C(\Omega\times [0,\tau])^3$ .

 Since $w_{n_k}\xrightharpoonup{k} 0$ weakly in ${\mathcal{F}}$,  we know $z_{n_k}\xrightarrow{k} 0$ by contradiction(cf. \cite{ref10}). That implies,
\begin{align}\label{eq-19}
&\int_0^t[v_{n_k}(x,r)-v(x,r)] dr \xrightarrow{k} 0\ \mathrm{uniformly\  in}\  \Omega\times [0,\tau].&
\end{align}

\vskip1mm Similarly,
\begin{align}
 &\|\nabla_x z_{n_k}(x,s)-\nabla_y z_{n_k}(y,t)\|\leq\left\|\int_0^s [\nabla_x w_{n_k}(x,r)-\nabla_y w_{n_k}(y,r)]dr\right\|+\left\|\int_s^t \nabla_y w_{n_k}(y,r)dr\right\| \nonumber &\\
 =&\left\|\int_0^s [\nabla_x v_{n_k}(x,r)-\nabla_y v_{n_k}(y,r)]-[\nabla_x v(x,r)-\nabla_y v(y,r)]dr\right\|
 +\left\|\int_s^t \nabla_y w_{n_k}(y,r)dr\right\| \nonumber &\\
 \leq &K\int_0^\tau \|Lv_{n_k}(\cdot,r)\|_{L^2(\Omega)}dr\|x-y\|^{\frac{1}{2}}
 +K\int_0^\tau \|Lv(\cdot,r)\|_{L^2(\Omega)}dr\|x-y\|^{\frac{1}{2}}+\left\|\int_s^t \nabla_y w_{n_k}(y,r)dr\right\| \nonumber &\\
 \leq &K\tau^{\frac{1}{2}}\left(\int_0^\tau \|Lv_{n_k}(\cdot,r)\|_{L^2(\Omega)}^2dr\right)^{\frac{1}{2}}\|x-y\|^{\frac{1}{2}}
 +K\tau^{\frac{1}{2}}\left(\int_0^\tau \|Lv(\cdot,r)\|_{L^2(\Omega)}^2dr\right)^{\frac{1}{2}}\|x-y\|^{\frac{1}{2}}+\left\|\int_s^t \nabla_y w_{n_k}(y,r)dr\right\| \nonumber & \\
 \leq &2KM^{\frac{1}{2}}\tau^{\frac{1}{2}}\|x-y\|^{\frac{1}{2}}+\left\|\int_s^t \nabla_y w_{n_k}(y,r)dr\right\|
  \leq 2KM^{\frac{1}{2}}\tau^{\frac{1}{2}}\|x-y\|^{\frac{1}{2}}+|t-s|^{\frac{1}{2}}K\|w_{n_k}(\cdot,r)\|_{\mathcal{F}} \nonumber & \\
  \leq &2KM^{\frac{1}{2}}\tau^{\frac{1}{2}}\|x-y\|^{\frac{1}{2}}+2MK|t-s|^{\frac{1}{2}}. \nonumber &
\end{align}

Hence, $\{\nabla_x z_{n_k}(x,s)\}$ is also relative compact in $C(\Omega\times [0,\tau])^3$.

Passing to a subsequence, we may assume that $\nabla_x z_{n_k}(x,s)\xrightarrow{k} G(x)$. Then $\nabla_x z_{n_k}(x,s)\xrightarrow{k} G(x)$ uniformly on $\Omega\times [0,\tau]$. On the other hand, (\ref{eq-19}) implies that $z_{n_k}(x,s)\xrightarrow{k} 0$ uniformly on $\Omega\times [0,\tau]$. Then we have $G(x)=\nabla_x 0=O_{3\times 3}$.
 Therefore,  $\nabla_x z_{n_k}(x,s)\xrightarrow{k} O_{3\times 3}$ uniformly in $\Omega\times [0,\tau]$, i.e., (\ref{newadd}) holds.

 Finally, we turn to showing that
\begin{align}\label{addmoi}
\Theta_{n_k}(s;0,x)\xrightarrow{k} \Theta(s;0,x).
\end{align}
\vskip1mm By Lemma \ref{th-jk6}, we obtain that
\begin{align}\label{eq-d15}
\Theta_{n_k}(s;0,x)-\Theta(s;0,x)&=\int_0^s \nabla_{\eta_{n_k}} v_{n_k}(\eta_{n_k}(s;0,x),s)\Theta_{n_k}(s;0,x)ds\nonumber&\\
&-\int_0^s\nabla_{\eta} v(\eta(s;0,x),s)\Theta(s;0,x)ds. &
\end{align}

 Then by Lemma \ref{lem-2}, we have
\begin{align}\label{eq-HY}
&\|\Theta_{n_k}(s;0,x)-\Theta(s;0,x)\|\nonumber&\\
&\leq {\int_0^s \|[\nabla_{\eta_{n_k}} v_{n_k}(\eta_{n_k}(s;0,x),s)-\nabla_{\eta} v_{n_k}(\eta(s;0,x),s)]\Theta_{n_k}(s;0,x)\|ds} \nonumber&\\
&+{\int_0^s\|\nabla_{\eta} v_{n_k}(\eta(s;0,x),s)[\Theta_{n_k}(s;0,x)-\Theta(s;0,x)]\|ds}\nonumber& \\
&+{\left\|\int_0^s  [\nabla_{\eta} v_{n_k}(\eta(s;0,x),s)-\nabla_{\eta} v(\eta(s;0,x),s)]\Theta(s;0,x)ds\right\|}\nonumber& \\
&=I_1+I_2+I_3.&
\end{align}

By Lemma \ref{lem-2}, we have
\begin{align}\label{eq-HYI1}
&I_1\leq \int_0^s K\|Lv_{n_k}(\cdot,s)\|_{L^2(\Omega)}\|\eta_{n_k}(s;0,x)-\eta(s;0,x)\|^{\frac{1}{2}}\|\Theta_{n_k}(s;0,x)\|ds \nonumber&\\
&\leq K\tilde{M}\tau^{\frac{1}{2}}\left(\int_0^s \|Lv_{n_k}(\cdot,s)\|_{L^2(\Omega)}^2ds\right)^{\frac{1}{2}}\|\eta_{n_k}(s;0,x)-\eta(s;0,x)\|^{\frac{1}{2}}_{C([0,\tau]:\overline{\Omega})} \nonumber&\\
&\leq K\tilde{M}\tau^{\frac{1}{2}}M^{\frac{1}{2}}\|\eta_{n_k}(s;0,x)-\eta(s;0,x)\|^{\frac{1}{2}}_{C([0,\tau]:\overline{\Omega})},
\end{align}
\begin{align}\label{eq-HYI2}
I_2&\leq \left\|\int_0^s\|\nabla_{\eta} v_{n_k}(\eta(s;0,x),s)\|\|\Theta_{n_k}(s;0,x)-\Theta(s;0,x)\| ds\right\| \nonumber& \\
&\leq \int_0^s\|\nabla_\eta v_{n_k}(\eta(s;0,x),s)\|\|\Theta_{n_k}(s;0,x)-\Theta(s;0,x)\| ds,
\end{align}

\vskip1mm On the other hand, we know $\|\eta_{n_k}(s;0,x)-\eta(s;0,x)\|_{C([0,\tau]:\overline{\Omega})}\xrightarrow{k} 0$ by {\bf(ii)}.  This leads to $I_1\rightarrow 0$ by (\ref{eq-HYI1}). By (\ref{newadd}), we know that $\nabla_x\left[\int_0^\tau v_{n_k}(x,r)-v(x,r)dr\right]\xrightarrow{k} O_{3\times 3}$ uniformly on $\Omega$. This implies that $I_3\rightarrow 0$.
\vskip1mm Hence, $I_1+I_3\rightarrow 0$ as $n_k \rightarrow \infty$. Then there exist $N=N(\varepsilon)$, such that $I_1+I_3<\varepsilon$ as $n_k>N$. Therefore, it follows from (\ref{eq-HY})  that
\begin{align}
\|\Theta_{n_k}(s;0,x)-\Theta(s;0,x)\|&\leq \varepsilon+ \int_0^s\|\nabla_\eta v_{n_k}(\eta(s;0,x),s)\|\|\Theta_{n_k}(s;0,x)-\Theta(s;0,x)\| ds. \nonumber
\end{align}
and Gronwall inequality[10, Lemma 1.1] implies that
\begin{align}
\|\Theta_{n_k}(s;0,x)-\Theta(s;0,x)\|&\leq \varepsilon e^{\int_0^s\|\nabla_\eta v_{n_k}(\eta(s;0,x),s)\| ds} \leq \varepsilon e^{K\tau^{\frac{1}{2}}(\int_0^\tau \|v_{n_k}(\cdot,s)\|^2_L ds)^{\frac{1}{2}}}\leq \varepsilon e^{K\tau^{\frac{1}{2}}M^{\frac{1}{2}}}.  & \nonumber
\end{align}
\vskip1mm So (\ref{addmoi}) is proved by letting $\varepsilon\rightarrow 0$.

\vskip1mm {\bf(iv)} For the sake of simplicity, in what follows, we denote $\Theta_{n_k}(s;0,x)$ and $\Theta(s;0,x)$ by $\Theta_{n_k}$ and $\Theta$, respectively. Then,

\begin{align}\label{tempory}
\|\Theta_{n_k}\Theta_{n_k}^T- \Theta\Theta^T\| &\leq \|\Theta_{n_k}\Theta_{n_k}^T- \Theta\Theta_{n_k}^T\|+\|\Theta\|\|\Theta_{n_k}^T-\Theta\Theta^T\| &\nonumber \\
&\leq \|\Theta_{n_k}- \Theta\|\|\Theta_{n_k}^T\|+\|\Theta\|\|\Theta_{n_k}^T-\Theta^T\| &\nonumber\\
&\leq \tilde{M}\|\Theta_{n_k}- \Theta\|+\tilde{M}\|\Theta_{n_k}^T-\Theta^T\|\xrightarrow{k} 0, &
\end{align}
since $\Theta_{n_k}\xrightarrow{k}\Theta$ and $\Theta_{n_k}^T\xrightarrow{k}\Theta^T$ by part {\bf (iii)}.

Similarly, there holds $\|\Theta_{n_k}^T\Theta_{n_k}- \Theta^T\Theta\|\xrightarrow{k} 0 $.

Now we prove that
\begin{align}\label{rotation}
R_{n_k}\xrightarrow{k} R
\end{align}
\vskip1mm Let $A_{n_k}=\Theta_{n_k}(s;0,x)\Theta_{n_k}^T(s;0,x)=\left(
                                                        \begin{array}{c}
                                                          a^{n_k}_{ij}(s;0,x) \\
                                                        \end{array}
                                                      \right)_{3\times 3},
                                                      A=\Theta(s;0,x)\Theta^T(s;0,x)=\left(
                                                        \begin{array}{c}
                                                          a_{ij}(s;0,x) \\
                                                        \end{array}
                                                      \right)_{3\times 3}
$.
\vskip1mm By (\ref{tempory}), we know $a^{n_k}_{ij}(s;0,x)\xrightarrow{k} a_{ij}(s;0,x)$ for $i,j=1,2,3$. Now, we simply denote  $a^{n_k}_{ij}(s;0,x)$ and  $a_{ij}(s;0,x)$ by $a^{n_k}_{ij}$ and $a_{ij}$, respectively.
\vskip1mm By Lemma \ref{th-jk6}, we obtain that
\begin{align}\label{eq-Fi2}
&{\rm det}(\Theta_{n_k}(\tau;0,x))=e^{-\int_0^\tau \sum\limits_{i=1}^3(v_{n_k}(\eta_{n_k}(s;t,x),s))_{i,x_i}ds}\neq 0. &
\end{align}
\vskip1mm By (\ref{eq-Fi2}), we know that $\Theta_{n_k}\Theta_{n_k}^T,\Theta\Theta^T\in SPD(3)$. Let $\lambda_{n_k}^{(1)}\geq \lambda_{n_k}^{(2)}\geq \lambda_{n_k}^{(3)}>0$, $\lambda^{(1)}\geq \lambda^{(2)}\geq \lambda^{(3)}>0$ be the eigenvalues of $\Theta_{n_k}\Theta_{n_k}^T,\Theta\Theta^T$ respectively.  By Lemma \ref{new}, we obtain $\lambda_{n_k}^{(i)}\rightarrow \lambda^{(i)} (i=1,2,3)$.
\vskip1mm By (\ref{eq-Fi2}), we obtain that
\begin{align}\label{cvb}
\det(A_{n_k})=&e^{-2\int_0^\tau \sum\limits_{i=1}^3(v_{n_k}(\eta_{n_k}(s;t,x),s))_{i,x_i}ds}\leq e^{2\int_0^\tau\|\nabla v(\cdot,s)\|ds}
\leq e^{2\int_0^\tau\|v_{n_k}(\cdot,s)\|_{[C^{1,\frac{1}{2}}(\Omega)]^3}ds}\nonumber \\
\leq & e^{2C\int_0^\tau\|v_{n_k}(\cdot,s)\|_{[H^3_0(\Omega)]^3}ds}
\leq e^{2K\int_0^\tau\|Lv_{n_k}(\cdot,s)\|_{L^2(\Omega)}ds}
\leq e^{2K\tau^{\frac{1}{2}}[\int_0^\tau\|Lv_{n_k}(\cdot,s)\|^2_{L^2(\Omega)}ds]^{\frac{1}{2}}}\nonumber \\
\leq & e^{2K\tau^{\frac{1}{2}}M^{\frac{1}{2}}}\doteq M_1<+\infty.&
\end{align}
\vskip1mm In a similar way, we can obtain that
\begin{align}\label{cvb1}
\det(A_{n_k})\geq \frac{1}{M_1}.
\end{align}
\vskip1mm By singularity decomposition theorem\cite{ref14}(See Appendix 2), we can find two orthogonal matrix $U_{n_k}$, $V_{n_k}$ such that
$\Theta_{n_k}=U_{n_k}S_{n_k}V_{n_k}^T$, where $S^{n_k}=diag(\sqrt{\lambda_{n_k}^{(1)}},\sqrt{\lambda_{n_k}^{(2)}},\sqrt{\lambda_{n_k}^{(3)}})$, $\lambda_{n_k}^{(1)},\lambda_{n_k}^{(2)},\lambda_{n_k}^{(3)}$ are eigenvalues of $\Theta_{n_k}\Theta_{n_k}^T$, $U_{n_k}$, $V_{n_k}$ are orthogonal eigenvectors of $\Theta_{n_k}\Theta_{n_k}^T$ and $\Theta_{n_k}^T\Theta_{n_k}$ respectively.
\vskip1mm Then, $A_{n_k}=\Theta_{n_k}\Theta_{n_k}^T=U_{n_k}(S_{n_k})^2U_{n_k}^T$ and $(A_{n_k})^{-1}=U_{n_k}(S_{n_k})^{-2}U_{n_k}^T$. Hence,
\begin{align}\label{ghjk}
\|(A_{n_k})^{-1}\|\leq&\|U_{n_k}\|\|(S_{n_k})^{-2}\|\|U_{n_k}^T\|
\leq \|U_{n_k}\|^2[\frac{1}{\lambda^{(1)}_{n_k}}+\frac{1}{\lambda^{(2)}_{n_k}}+\frac{1}{\lambda^{(3)}_{n_k}}]&\nonumber \\
\leq&\|U_{n_k}\|^2\frac{\lambda^{(1)}_{n_k}\lambda^{(2)}_{n_k}+\lambda^{(1)}_{n_k}\lambda^{(3)}_{n_k}+\lambda^{(2)}_{n_k}\lambda^{(3)}_{n_k}}{\lambda^{(1)}_{n_k}\lambda^{(2)}_{n_k}\lambda^{(3)}_{n_k}}
=\|U_{n_k}\|^2\frac{\lambda^{(1)}_{n_k}\lambda^{(2)}_{n_k}+\lambda^{(1)}_{n_k}\lambda^{(3)}_{n_k}+\lambda^{(2)}_{n_k}\lambda^{(3)}_{n_k}}{\det(A_{n_k})}&\nonumber\\
\leq&\|U_{n_k}\|^2\frac{[\lambda^{(1)}_{n_k}+\lambda^{(2)}_{n_k}+\lambda^{(3)}_{n_k}]^2}{\det(A_{n_k})}
=\|U_{n_k}\|^2\frac{[tr(A_{n_k})]^2}{\det(A_{n_k})}
\leq \|U_{n_k}\|^2\frac{\|A_{n_k}\|^2}{\det(A_{n_k})}&\nonumber \\
\leq&27\tilde{M}^2M_1\doteq M_2<+\infty,&
\end{align}
by (\ref{cvb}), (\ref{cvb1}) and (\ref{cvb3}),where $tr(A)$ denote the trace of matrix $A$.
\vskip1mm Similarly , we know that $\|A^{-1}\|\leq M_2$.
\vskip1mm By (\ref{ghjk}), we obtain that
\begin{align}\label{nerf}
\|(A_{n_k})^{-1}-A^{-1}\|=&\|A_{n_k}^{-1}(A-A_{n_k})A^{-1}\|\leq \|A_{n_k}^{-1}\|\|A-A_{n_k}\|\|A^{-1}\|\xrightarrow{k} 0,&
\end{align}
since $A-A_{n_k}=\Theta_{n_k}\Theta_{n_k}^T-\Theta\Theta^T\rightarrow 0$ by (\ref{tempory}).

\vskip1mm Hence, $A_{n_k}^{-1}\xrightarrow{k} A^{-1}$.
\vskip1mm Since $A_{n_k}^{-1}=U_{n_k}(S_{n_k})^{-2}U_{n_k}^T=\left[U_{n_k}S_{n_k}^{-1}U_{n_k}^T\right]\left[U_{n_k}S_{n_k}^{-1}U_{n_k}^T\right]\triangleq B_{n_k}B_{n_k}$, $A^{-1}=US^{-2}U^T=\left[US^{-1}U^T\right]\left[US^{-1}U^T\right]\triangleq BB$, then $B_{n_k}$ and $B$ are positive definite symmetric matrixes. By Lemma \ref{new}, we know $\det(B_{n_k})=\frac{1}{\sqrt{\lambda_{n_k}^{(1)}\lambda_{n_k}^{(2)}\lambda_{n_k}^{(3)}}}\rightarrow \frac{1}{\sqrt{\lambda^{(1)}\lambda^{(2)}\lambda^{(3)}}}=\det(B)$. Since $B_{n_k}$, $B\in SPD(3)$, the Minkowskii inequality implies
\begin{align}
\left[\det(B_{n_k}+B)\right]^{\frac{1}{3}}\geq \left[\det(B_{n_k})\right]^{\frac{1}{3}}+\left[\det(B)\right]^{\frac{1}{3}}\geq \left[\det(B)\right]^{\frac{1}{3}}.
\end{align}

This is, $\det(B_{n_k}+B)\geq \det(B)=\frac{1}{\lambda^{(1)}\lambda^{(2)}\lambda^{(3)}}>0$.
\vskip1mm Further more, we have
\begin{align}\label{nerfhj1}
\|B_{n_k}+B\|\leq &\|B_{n_k}\|+\|B\|
\leq\|U_{n_k}\|^2[\frac{1}{\sqrt{\lambda_{n_k}^{(1)}}}+\frac{1}{\sqrt{\lambda_{n_k}^{(2)}}}+\frac{1}{\sqrt{\lambda_{n_k}^{(3)}}}]+\|U\|^2[\frac{1}{\sqrt{\lambda^{(1)}}}+\frac{1}{\sqrt{\lambda^{(2)}}}+\frac{1}{\sqrt{\lambda^{(3)}}}]&\nonumber\\
=&\|U_{n_k}\|^2\frac{\sqrt{\lambda_{n_k}^{(1)}\lambda_{n_k}^{(2)}}+\sqrt{\lambda_{n_k}^{(1)}\lambda_{n_k}^{(3)}}+\sqrt{\lambda_{n_k}^{(2)}\lambda_{n_k}^{(3)}}}{\sqrt{\lambda_{n_k}^{(1)}\lambda_{n_k}^{(2)}\lambda_{n_k}^{(3)}}}
+\|U\|^2\frac{\sqrt{\lambda^{(1)}\lambda^{(2)}}+\sqrt{\lambda^{(1)}\lambda^{(3)}}+\sqrt{\lambda^{(2)}\lambda^{(3)}}}{\sqrt{\lambda^{(1)}\lambda^{(2)}\lambda^{(3)}}}&\nonumber\\
\leq&\|U_{n_k}\|^2\frac{\lambda_{n_k}^{(1)}+\lambda_{n_k}^{(2)}+\lambda_{n_k}^{(3)}}{\sqrt{\lambda_{n_k}^{(1)}\lambda_{n_k}^{(2)}\lambda_{n_k}^{(3)}}}+\|U\|^2\frac{\lambda^{(1)}+\lambda^{(2)}+\lambda^{(3)}}{\sqrt{\lambda^{(1)}\lambda^{(2)}\lambda^{(3)}}}
=\|U_{n_k}\|^2\frac{tr(A_{n_k})}{\sqrt{\det(A_{n_k})}}+\|U\|^2\frac{tr(A)}{\sqrt{\det(A)}}&\nonumber \\
\leq& 9M_2\sqrt{M_1}+9M_2\sqrt{M_1}\leq 18M_2\sqrt{M_1}\doteq M_3<+\infty.&
\end{align}
\vskip1mm By (\ref{nerf}), we obtain that
\begin{align}\label{nerf1}
A_{n_k}^{-1}-A^{-1}=B_{n_k}^2-B^2=(B_{n_k}+B)(B_{n_k}-B)\xrightarrow{k} 0.
\end{align}
\vskip1mm By (\ref{nerf1}) and Lemma \ref{newtg}, we obtain that $(\Theta_{n_k}\Theta_{n_k}^T)^{-\frac{1}{2}}=B_{n_k}\rightarrow B=(\Theta\Theta^{T})^{-\frac{1}{2}}$.
\vskip1mm Based on the above calculation, we obtain that
\begin{align}\label{nerf2}
&\|\Theta_{n_k}(\Theta_{n_k}\Theta_{n_k}^T)^{-\frac{1}{2}}-\Theta(\Theta\Theta)^{-\frac{1}{2}}\|&\nonumber \\
\leq&\|\Theta_{n_k}\|\|(\Theta_{n_k}\Theta_{n_k}^T)^{-\frac{1}{2}}-(\Theta\Theta)^{-\frac{1}{2}}\|+\|\Theta_{n_k}-\Theta\|\|(\Theta\Theta)^{-\frac{1}{2}}\|\xrightarrow{k} 0.&
\end{align}
\vskip1mm So, $R_{n_k}\xrightarrow{k} R$.\qed
\section{Proof of Theorem\ref{th-A10}}
The aim of this section is to prove Theorem\ref{th-A10}. For this purpose, we let
\begin{align}
\inf_{v\in{\mathcal{F}}} H(v)=\bar{H}.
\end{align}

\noindent{\bf Proof of Theorem \ref{th-A10}:}\ \ Let $\{v_n\}\subset \mathcal{F}$ be a minimizing sequence of $H(v)$, that is,
\begin{align}\label{conclu}
\lim_{n\rightarrow\infty} H(v_n)=\bar{H}.
\end{align}

Then, there exists a constat $M>0$ such that $\|v_n\|_{\mathcal{F}}^2\leq M$, otherwise, there is a contradiction by $\mathop {\lim }\limits_{\|v\|_{\mathcal{F}}\rightarrow\infty}H(v)=+\infty$.

 It follows from Lemma \ref{pd7}.{\bf(i)} that there exists a weakly convergent subsequence $\{v_{n_k}\}$ of $\{v_n\}$ such that
 \begin{align}
 v_{n_k}\xrightharpoonup{k} \hat{v}, \ \ \mathrm{for \ some}\ \hat{v}\in \mathcal{F}.
 \end{align}

Now, we claim that
\begin{align}\label{eq-cD23}
&\lim_{n_k\rightarrow \infty}\inf {H}(v_{n_k})\geq {H}(\hat{v}).&
\end{align}

By Lemma \ref{pd7}{\bf(i)}, we know that
\begin{align}\label{eq-dhk10}
&M\geq\lim_{n_k\rightarrow \infty} \inf \int_0^\tau \|Lv_{n_k}(\cdot,t)\|^2_{L^2(\Omega)}dt\geq \int_0^\tau \|L\hat{v}(\cdot,t)\|^2_{L^2(\Omega)}dt, &
\end{align}

Further more, by Lemma \ref{pd7} {\bf(ii)}, {\bf(iii)}, {\bf(iv)} and definition of $h_{n_k}(x)$ and $h(x)$ , we have
\begin{align}\label{eq-derh13}
&h_{n_k}(x)=\eta_{n_k}(0;\tau,x)\xrightarrow{k} \eta(0;\tau,x)=h(x) \ \ {\rm for} \ \ \forall x \in \Omega ,&
\end{align}
and
\begin{align}\label{acc}
R_{n_k}\xrightarrow{k} R.
\end{align}


\vskip1mm It follows from (\ref{eq-derh13}) and (\ref{acc}) that
\begin{align}
\|T\diamond h_{n_k}(x)-T\diamond h(x)\|&=\|R_{n_k} [T\circ h_{n_k}(x)] R_{n_k}^T-R [T\circ h(x)] R^{ T}\| \nonumber &\\
&\leq\|R_{n_k} [T\circ h_{n_k}(x)] R_{n_k}^T-R [T\circ h_{n_k}(x)] R_{n_k}^T\| \nonumber &\\
&+\|R [T\circ h_{n_k}(x)] R_{n_k}^T-R [T\circ h(x)] R_{n_k}^T\| \nonumber &\\
&+\|R [T\circ h(x)] R_{n_k}^T-R [T\circ h(x)] R^{ T}\| \nonumber &\\
&\leq\|R_{n_k}-R\|\| [T\circ h_{n_k}(x)] R_{n_k}^T\| \nonumber &\\
&+\|R\|\| [T\circ h_{n_k}(x)] - [T\circ h(x)]\|\|R_{n_k}^T\| \nonumber &\\
&+\|R [T\circ h(x)]\|\| R_{n_k}^T- R^{T}\|\xrightarrow{k} 0 \ \mathrm{on} \ x\in \Omega\setminus h^{-1}(\Delta_T). \nonumber &
\end{align}
\vskip1mm Hence,  $\|T\diamond h_{n_k}(\cdot)-D(\cdot)\|^2\xrightarrow{k} \|T\diamond h(\cdot)-D(\cdot)\|^2$ on $x\in \Omega\setminus h^{-1}(\Delta_T)$.

By Remark \ref{newremark},
we know that $h:\Omega\rightarrow\Omega$ is a $1$-to-$1$ and onto, which ensure the existence of $h^{-1}(x)$. Therefore, $h^{-1}(\Delta_T)$ is a set of Lebesgue measure zero.
\vskip1mm On the other hand, we know
\begin{align}
\|T\diamond h_{n_k}(\cdot)-D(\cdot)\|^2\leq \mathop {\max }\limits_{x\in \Omega} \|T(x)-D(x)\|^2=J\in L^1(\Omega).
\end{align}

By the Lebesgue Dominant Convergence Theorem, we obtain that
\begin{align}\label{eq-D22}
&\|T\diamond h_{n_k}(\cdot)-D(\cdot)\|_{L^2(\Omega)}^2\xrightarrow{k} \|T\diamond h(\cdot)-D(\cdot)\|_{L^2(\Omega)}^2. &
\end{align}
\vskip1mm Therefore, (\ref{eq-D22}) and (\ref{eq-dhk10}) implies that (\ref{eq-cD23}) holds.

\vskip1mm Since $\hat{v}\in {\mathcal{F}}$,
\begin{align}\label{ghkm1}
H(\hat{v})\geq \inf_{v\in{\mathcal{F}}} H(v)=\bar{H}.
\end{align}

Combining (\ref{conclu}) and (\ref{eq-cD23}), we see that
\begin{align}\label{ghkm2}
&\bar{H}=\lim_{n_k\rightarrow \infty}\inf {H}(v_{n_k})\geq H(\hat{v})\geq\bar{H}.&
\end{align}

That is,
\begin{align}\label{ghkm3}
&H(\hat{v})=\bar{H}=\inf_{v\in{\mathcal{F}}} H(v).&
\end{align}
\vskip1mm For the above minimizer $\hat{v}\in\mathcal{F}$, by Lemma \ref{th-5}, we know that there exists an unique $\hat{\eta}(s;t,x)$ such that
\begin{align}
\frac{d\hat{\eta}(s;t,x)}{ds}=\hat{v}(\hat{\eta}(s;t,x),s), \ \ \ \ \hat{\eta}(t;t,x)=x.
\end{align}
\vskip1mm Hence, the mapping $\hat{h}(x)=\hat{\eta}(0;\tau,x):\Omega\rightarrow\Omega$ is what we want. Moreover, by {\bf(ii)} in Lemma \ref{pd7}, we know $\hat{h}(x)\in [C^{1,\frac{1}{2}}(\Omega)]^3$ with the derivative satisfies (\ref{eq-5}).\qed

\vskip1mm  {\bf Acknowledgements.}{\ \ This work was supported by NSFC under grant 11471331 and National Center for Mathematics and Interdisciplinary Sciences.}\\

{\bf\large  Appendix: The calculation of $(\mathbf{A}\mathbf{A}^T)^{-\frac{1}{2}}$}
\begin{appendix}

\vskip1mm By singularity decomposition theorem\cite{ref14}, for any matrix $\mathbf{A}\in \mathbb{R}^{n\times n}$, it could be broken down into the product of three matrixes:
\begin{align}\label{eq-q}
&\mathbf{A}=\mathbf{U}\mathbf{S}\mathbf{V}^T,&
\end{align}
where $\mathbf{U}^T\mathbf{U}=\mathbf{U}\mathbf{U}^T=\mathbf{I},\mathbf{V}^T\mathbf{V}=\mathbf{V}\mathbf{V}^T=\mathbf{I}$; The columns of $\mathbf{U}$ are orthogonal eigenvectors[7] of $\mathbf{A}\mathbf{A}^T$, the columns of $\mathbf{V}$ are orthogonal eigenvectors of $\mathbf{A}^T\mathbf{A}$. $\mathbf{S}={\rm diag}(\sigma^{(1)},\sigma^{(2)},\cdots,\sigma^{(n)})$, and $\sigma^{(i)}=\sqrt{\lambda^{(i)}}(i=1,2,\cdots,n)$. Where ${\lambda^{(i)}}(i=1,2,\cdots,n)$ are eigenvalues of $\mathbf{A}^T\mathbf{A}$. Here we call $\mathbf{U},\mathbf{V}$ orthogonal matrixes of $\mathbf{A}\mathbf{A}^T$ and $\mathbf{A}^T\mathbf{A}$.

\vskip1mm Based on the theory mentioned above, we focus on the calculation of $(\mathbf{A}\mathbf{A}^T)^{-\frac{1}{2}}$.
\vskip1mm By (\ref{eq-q}), we obtain
\begin{align}\label{eq-r}
\mathbf{A}\mathbf{A}^T&=\mathbf{U}\mathbf{S}\mathbf{V}^T\mathbf{V}\mathbf{S}\mathbf{U}^T=\mathbf{U}\mathbf{S}^2\mathbf{U}^T
=\mathbf{U}\mathbf{S}\mathbf{U}^T\mathbf{U}\mathbf{S}\mathbf{U}^T.&
\end{align}
\vskip1mm Therefore, by the uniqueness of square root of matrix in [14], we obtain that
\begin{align}\label{eq-e}
(\mathbf{A}\mathbf{A}^T)^{\frac{1}{2}}&=\mathbf{U}\mathbf{S}\mathbf{U}^T.&
\end{align}
\vskip1mm  By (\ref{eq-e}), we know
\begin{align}\label{eq-f}
(\mathbf{A}\mathbf{A}^T)^{-\frac{1}{2}}&=\mathbf{U}\mathbf{S}^{-1}\mathbf{U}^T.&
\end{align}

\end{appendix}

\end{document}